\documentstyle{article}
\input amssym.def
\input amssym.tex

\newcommand{\forces}{\Vdash}

\newcommand{\lbv}{[\![} 
\newcommand{\rbv}{]\!]}
\newcommand{\V}{{\bf V}} 
\newcommand{\lesdot}{\mathrel{\mathord{<}\!\!\raise 0.8
pt\hbox{$\scriptstyle\circ$}}}  
 
\newcommand{\non}{{\rm {\bf non}}\/} 
\newcommand{\con}{{\frak c}} 
\newcommand{\dominating}{{\frak d}} 
\newcommand{\unbounded}{{\frak b}}
 
\newcommand{\can}{2^{\textstyle \omega}} 
\newcommand{\fs}{2^{\textstyle <\!\omega}} 
\newcommand{\baire}{\omega^{\textstyle \omega}} 
\newcommand{\iso}{[\omega]^{\textstyle \omega}} 
\newcommand{\fsuo}{[\omega]^{\textstyle <\!\omega}} 
\newcommand{\fseo}{\omega^{\textstyle <\!\omega}} 
\newcommand{\lh}{{\rm lh}\/} 
\newcommand{\rest}{{\restriction}}
\newcommand{\lra}{\longrightarrow}
 
\newcommand{\suc}{{\rm succ}} 
\newcommand{\dom}{{\rm dom}} 
\newcommand{\rng}{{\rm rng}}
\newcommand{\nor}{{\rm {\bf nor}}\/} 

\newcommand{\random}{{\Bbb B}}

\newcommand{\K}{{\cal K}}
\newcommand{\k}{{\Bbb K}}

\renewcommand{\P}{{\cal P}}
\newcommand{\p}{{\Bbb P}}

\newcommand{\q}{{\Bbb Q}}

\renewcommand{\t}{{\Bbb T}}   
\newcommand{\X}{{\cal X}}
\newcommand{\pk}{{\cal P}_k}
\newcommand{\icso}{[\omega]^{\textstyle *}}
\newcommand{\rall}{R^\forall_k}
\newcommand{\rexists}{R^\exists_k}
\newcommand{\cont}{{\rm cont}}
\newcommand{\iseo}{\omega^{\textstyle{<}\omega}}
\newcommand{\uh}{{{\rm uh}}}

\newcommand{\QED}{\hfill\hspace{0.2in}\vrule width 6pt height 6pt depth 0pt 
\vspace{0.1in}} 
\newcommand{\Proof}{\noindent{\sc Proof} \hspace{0.2in}} 
\newtheorem{theorem}{Theorem}[section] 
 
\newtheorem{claim}{Claim}[theorem]
\newtheorem{mainclaim}{Main Claim}[theorem]
\newtheorem{notation}[theorem]{Basic Definitions}
\newtheorem{proposition}[theorem]{Proposition} 
\newtheorem{corollary}[theorem]{Corollary} 
\newtheorem{problem}[theorem]{Problem} 
\newtheorem{definition}[theorem]{Definition}
\newtheorem{example}[theorem]{Example}
\title{Localizations of infinite subsets of $\omega$} 
\author{
{\bf Andrzej Ros{\l}anowski} {\thanks{The research was partially
supported by Polish Committee of Scientific Research, Grant KBN 654/2/91}}\\
Institute of Mathematics\\
The Hebrew University of Jerusalem\\
Jerusalem, Israel\\
and\\
Mathematical Institute of Wroc{\l}aw University\\
50384 Wroc{\l}aw, Poland
\and
{\bf Saharon Shelah}\thanks{The second author would like to thank
{\em Basic Research Foundation} of The Israel Academy of
Sciences and Humanities for partial support. Publication number 501.}\\
Institute of Mathematics\\
The Hebrew University of Jerusalem\\
Jerusalem, Israel\\
and\\
Department of Mathematics\\
Rutgers University\\
New Brunswick, NJ 08854, USA
}
 
\date{done: September 1992\\
printed: \today}
\begin{document} 
\maketitle 

\eject

\setcounter{section}{-1}
\section{Introduction}

\subsection{Preliminaries}
When we say {\em reals} we mean one of the nicely defined Polish
spaces or their (finite or countable) products like: the real line ${\Bbb R}$,
the Cantor space $\can$, the Baire space $\baire$ or the space of infinite
sets of integers $\iso$.\\ 
In the present paper we are interested in properties of forcing notions (or,
generally, extensions of models of ZFC) which measure in a sense the distance
between the ground model reals and the reals in the extension. In particular
we look at the ways the ``new'' reals can be localized (or: aproximated)
by ``old'' reals. There are two extreme cases here: {\em there are no new
reals} and {\em the old reals are countable}. However, between these two
extremes we have a wide spectrum of properties among which the localizations
by slaloms seem to be the most popular. A systematic study of slaloms and
related localization properties and cardinal invariants was presented in
\cite{Bar1}. 
\smallskip

\noindent{\em A slalom} is a function $S:\omega\lra\fsuo$ such that
$(\forall n\!\in\!\omega)(|S(n)|=n+1)$. We say that a slalom $S$ {\em
localizes} a function $f\in\baire$ whenever $(\forall n\!\in\!\omega)(f(n)\in
S(n))$. In this situation we can think that the slalom $S$ is an approximation
of the function $f$. It does not determine the function but it provides some
bounds on possible values of $f$.  Bartoszy\'nski, Cicho\'n, Kamburelis et al.
studied the localization by slaloms and those investigations gave the
following surprising result.
\begin{theorem}[Bartoszy\'nski, \cite{Bar2}]
\label{bartoszynski}
\hspace{0.1in} Suppose that $\V\subseteq\V'$ are models of ZFC.  Then the
following conditions are equivalent: 
\begin{enumerate}
\item Any function from $\baire\cap\V'$ can be localized by a
slalom from \V.
\item Any Borel (Lebesgue) null set coded in $\V'$ can be covered by a
Borel null set coded in \V. \QED
\end{enumerate}
\end{theorem}
On localizations by slaloms see Chapter VI of \cite{Sh:b} too; other
localizations of slalom-like type appeared in \cite{GoSh:448}.

\noindent A stronger localization property was considered in \cite{NeRo}. Fix
a natural number $k\geq 2$. By a {\em $k$-tree} on $\omega$ we mean a tree
$T\subseteq\iseo$ such that each finite node in $T$ has at most $k$ immediate
successors in $T$. We say that a $k$-tree $T$ {\em localizes} a function
$f\in\baire$ whenever $f$ is a branch through $T$ (i.e.  $(\forall
n\!\in\!\omega)(f\rest n\in T)$). Clearly, if any function from
$\baire\cap\V'$ can be localized by a $k$-tree from $\V$, $k\geq 2$,
$\V\subseteq\V'$ then each function from $\baire\cap\V'$ can be localized by a
slalom from $\V$.  Moreover the localization by a $k$-tree implies the
localization by a $k+1$-tree (but not conversely).

\noindent In the definition below we formulate general localization properties
for Polish spaces  $X,Y$. In practice, however, these spaces will be various
examples of reals only. 
\begin{definition}
\label{defloc}
Assume that $X,Y$ are Polish spaces and $R\subseteq X\times Y$ is a Borel
relation. Suppose that $\V\subseteq\V'$ are models of ZFC and that all
parameters we need are in $\V$. We say that the pair $(\V,\V')$ has {\em the
property of the $R$-localization} if  
\[(\forall x\!\in\! X\!\cap\!\V')(\exists y\!\in\! Y\!\cap\!\V)
((x,y)\in R)\] 
If $x\in X\cap\V'$, $y\in Y\cap\V$ and $(x,y)\in R$ they we say that {\em $y$
$R$-localizes $x$}. 
\end{definition}
In the examples we gave earlier $X$ was $\baire$ and $Y$ was the space of
slaloms or the space of all $k$-trees, respectively. The respective relations
should be obvious. Those localizations were to approximate functions in an
extension by objects from the ground model. They are not useful if we
consider infinite subsets of $\omega$. Though each member of $\iso$ can be
identified with its increasing enumeration, the localization (either by
slaloms or by $k$-trees) of the enumeration does not provide satisfactory
information on successive points of the set. The localization gives us
``candidates'' for the $n$-th point of the set but the same candidates can
appear several times for distinct $n$.  That led to a suggestion that we
should consider disjoint subsets of $\omega$ as sets of ``candidates'' for
successive points of the localized set (the approach was suggested by B.
W\c eglorz). Now we have two possibilities. Either we can demand that each set
from the localization contains a limited number of members of the localized
set or we can postulate that each intersection of that kind is large.
Localizations of this kind are studied in section 1. In the second section we
investigate localizations of infinite subsets of $\omega$ by sets of integers
from the ground model. These localizations might be thought as localizations
by partitions of $\omega$ into successive intervals. A starting point for our
considerations was the following observation. 
\begin{proposition}
Suppose that $\V\subseteq\V'$ are (transitive) models of ZFC. Then:
\begin{enumerate}
\item $\V\cap\baire$ is unbounded in $\V'\cap\baire$ if and only if

for every set $X\in\iso\cap\V'$ there exists a set $Yin\iso\cap\V$ such that
{\bf infinitely often} between two successive points of $Y$ there are at least
2 points of $X$.
\item $\V\cap\baire$ is dominating in $\V'\cap\baire$ if and only if

for every set $X\in\iso\cap\V'$ there exists a set $Y\in\iso\cap\V$ such that
{\bf for all but finitely many} pairs of two successive points of $Y$ there
are at least 2 points of $X$ between them.
\end{enumerate}
\end{proposition}
Now we try to replace the quantifier {\em for infinitely many} above by
stronger quantifiers (but still weaker than {\em for all but finitely many}),
like {\em for infinitely many $n$, for both $n$ and $n+1$}.
Finally, in section 3 we formulate several corollaries to the results of
previous sections for cardinal invariants related to the notions we study.

\subsection{Notation} Our notation is rather standard and
essentially compatible with that of \cite{Jec} and \cite{Kun}. In forcing
considerations, however, we will use the convention that {\em a stronger
condition is the greater one}.
\begin{notation}
\begin{enumerate}
\item A tree on $\omega$ is a set $T\subseteq\fseo$ closed under initial
segments. For the tree $T$ the body $[T]$ of $T$ is the set
\[\{x\in\baire:(\forall l\in\omega)(x\rest l\in T)\}.\] 
If $t\in T$ then $\suc_T(t)=\{s\in T: t\subseteq s\ \&\ \lh(t)+1=\lh(s)\}$.
\item By a model of ZFC we will mean {\em a transitive model} of (enough
of) ZFC. Models of ZFC will be denoted by $\V$, $\V'$ etc. 
\item We will be interested in extensions of models, i.e. in pairs $(\V,\V')$
of models such that $\V\subseteq\V'$. If a property of an
extension is defined then we extend this definition to notions of forcing.
We say that a notion of forcing $\p$ has the property whenever for any
generic filter $G\subseteq\p$ over $\V$ the extension $\V\subseteq\V[G]$ has
the considered property.
\item We will use the quantifiers $(\forall^\infty n)$ and $(\exists^\infty
n)$ as abbreviations for 
\[(\exists m\in\omega)(\forall n>m)\quad\mbox{ and }\quad(\forall m\in\omega)
(\exists n>m),\] 
respectively.  
\item The Baire space $\baire$ of all functions from $\omega$ to
$\omega$ is endowed with the partial order $\leq^*$:
\[f\leq^* g\ \iff\ (\forall^\infty n)(f(n)\leq g(n)).\]
A family $F\subseteq\baire$ is unbounded in $(\baire,\leq^*)$ if
\[\neg(\exists g\in\baire)(\forall f\in F)(f\leq^* g)\] 
and it is dominating in $(\baire,\leq^*)$ if 
\[(\forall g\in\baire)(\exists f\in F)(g\leq^* f).\] 
\item The unbounded number $\unbounded$ is the minimal size of an unbounded
family in $(\baire,\leq^*)$; the dominating number $\dominating$ is the
minimal size of a dominating family in that order.
\item The size of the continuum is denoted by $\con$, $\icso$ stands for the
family of infinite co-infinite subsets of $\omega$. 
\end{enumerate}
\end{notation}

\subsection{Acknowledgement} We would like to thank Professor Uri Abraham
for his helpful comments on the paper.

\section{$\rall$, $\rexists$ - localizations}

In this section we show that a localization of infinite subsets of $\omega$
suggested by B. W\c eglorz implies that the considered extension adds no new
real.
\begin{definition}
\label{defrel}
\begin{enumerate}
\item A partition of $\omega$ into finite sets is a sequence $\langle
K_n:n\in\omega\rangle$ of disjoint finite sets such that
$\bigcup_{n\in\omega}K_n=\omega$. 
\item $\pk$ is the set of all partitions $\langle K_n:n\in\omega\rangle$ of
$\omega$ into finite sets such that $(\forall n\!\in\!\omega)(|K_n|>k)$.

{\em [Note that $\pk$ is a $\Pi^0_2$-subset of $(\fsuo)^{\textstyle \omega}$
so it is a Polish space.]}
\item We define relations $\rall\subseteq\icso\times\pk$ and
$\rexists\subseteq\icso\times\P_{k+1}$ by 
\[(X,\langle K_n:n\!\in\!\omega\rangle)\in\rall\ \equiv\
(\forall^\infty n\!\in\!\omega)(|X\cap K_n|\leq k)\]
\[(X,\langle K_n:n\!\in\!\omega\rangle)\in\rexists\ \equiv\
(\exists^\infty n\!\in\!\omega)(|X\cap K_n|\leq k).\]
Their complements (in $\icso\times\pk$, $\icso\times\P_{k+1}$) are denoted
by $c\rall$, $c\rexists$, respectively. 
\end{enumerate}
\end{definition}
If we want to approximate an infinite co-infinite subset of $\omega$ by an
object in a given model we can look for a separation of distinct members of
the set by a sequence of sets from the model. Thus we could ask if it is
possible to find a partition of $\omega$ (in $\V$) such that the localized set
is {\em a partial selector} of the partition. More generally we may ask for
$\rall$--localization; recall definition \ref{defloc}. Thus the 
$\rall$-localization property means that for every infinite set of integers
$X$ from the extension there exists a partition $\langle K_n:
n\in\omega\rangle\in\pk$ from the ground model such that for almost all
$n\in\omega$ the intersection $X\cap K_n$ is of size at most $k$. The following
result shows that the $\rall$-localization fails if we add new reals.

\begin{theorem}
\label{trivial}
\hspace{0.1in} Suppose that $\V\subseteq\V'$ are models of ZFC such that
$\V\cap\can\neq\V'\cap\can$. Then there is a set $X\in\icso\cap\V'$
such that for {\bf no} $k\in\omega$ there is a partition $\langle
K_n:n\in\omega\rangle\in\V$ of $\omega$ such that
\[(\forall n\!\in\!\omega)(|X\cap K_n|\leq k \ \&\ \ |K_n|>k).\]
Consequently the extension $\V\subseteq\V'$ does not have the
$\rall$-localization property (for any $k$).
\end{theorem}

\Proof Let $x\in\can$ be a new real (i.e. $x\in\V'\setminus\V$) and let
\[X=\fseo\setminus\{x\rest i:i\in\omega\}.\] 
As we can identify $\omega$ with $\fseo$ we may think that $X\in\icso$. Now
the following claim works.

\begin{claim}
\label{lemat}
\hspace{0.1in} Suppose that $x\in\can$, $x\notin\V$. Let
$\langle K_n:n\in\omega\rangle\in\V$ be a sequence of disjoint subsets of
$\fseo$ such that $(\forall n\!\in\!\omega)(|K_n|>k)$, $k\in\omega$. Then
for some $n\in\omega$ the set $K_n\setminus\{x\rest i:i\in\omega\}$ has at
least $k+1$ points.
\end{claim}

\noindent{\em Proof of Claim:}\ \ \ Suppose not. Thus for each
$n\in\omega$ we have 
\[|K_n\setminus\{x\rest i: i\in\omega\}|\leq k\quad\mbox{ and thus}\quad
K_n\cap\{x\rest i: i\in\omega\}\neq\emptyset.\]  
First note that each $K_n$ is finite. If not then the tree 
\[\{s\in\fs:(\exists t\in K_n)(s\subseteq t)\}\]
has exactly one infinite branch - the branch is $x$. As the tree is in $\V$ we
would get $x\in\V$. 

\noindent Let
\[u(n)=\max\{\lh(s): s\in K_n\}\quad\mbox{ and }\quad d(n)=\min\{\lh(s): s\in
K_n\}\] 
(remember that each $K_n$ is finite). Choose an increasing sequence $\langle
n_\ell:\ell{<}\omega\rangle$ (in $\V$) such that
\[(\forall\ell\in\omega)(u(n_\ell)<d(n_{\ell+1}))\]
(possible as the $K_n$'s are disjoint). Let
\[F_\ell=\{s\rest d(n_\ell): s\in K_{n_\ell}\}.\]
Note that $x\rest d(n_\ell)\in F_\ell$ as an initial segment of $x$ belongs to
$K_{n_\ell}$. Moreover $|F_\ell|\leq k+1$ as only one element of $F_\ell$ may
be an initial segment of $x$ and above each member of $F_\ell$ there is an
element of $K_{n_\ell}$ (remember $| K_{n_\ell}\setminus\{x\rest i:
i\in\omega\}|\leq k$). Clearly $\langle F_\ell: \ell\in\omega\rangle\in\V$.
Consider the set
\[A=\{y\in\can: (\forall\ell\in\omega)(y\rest d(n_\ell)\in F_\ell)\}.\]
It is a finite set from $\V$. But $x\in A$ -- a contradiction. \QED

\noindent Thus the $\rall$-localization is the trivial one. The complementary
$c\rall$-localization is not of special interest either. Every extension
$\V\subseteq\V'$ has the $cR^{\forall}_0$-localization. The description of
the $c\rall$-localization for $k>0$ is given by the following observation.

\begin{proposition}
\label{unbounded}
\hspace{0.1in}Let $\V\subseteq\V'$ be an extension of models of
ZFC. Then the following conditions are equivalent:
\begin{enumerate}
\item For each $k>0$ the pair $(\V,\V')$ has the $c\rall$-localization.
\item For some $k>0$ the pair $(\V,\V')$ has the $c\rall$-localization.
\item $\V\cap\baire$ is unbounded in $\V'\cap\baire$.
\end{enumerate}
\end{proposition}

\Proof $2\Rightarrow 3$\ \ \ Let an increasing function
$f\in\baire\cap\V'$ be given. Take an increasing function
$f'\in\V'\cap\baire$ such that 
\[(\forall n\in\omega)(f(f'(n))+1<f'(n+1)).\]
Consider the set $\rng(f')\in\icso$.  Let $\langle K_n: n\in\omega\rangle\in
\V$ be the partition of $\omega$ given for this set by the
$c\rall$-localization.  Let $k_n=\min K_n$ and put $g(k_n)=1+\max K_n$. 
Extend $g$ to $\omega$ putting $g(m)=0$ if $m\notin\{k_n:n\in\omega\}$.
Clearly $g\in\V$.

\noindent Note that $|K_n\cap\rng(f')|>k>0$ implies that $f(k_n)<g(k_n)$.
Consequently 
\[(\exists^\infty m\!\in\!\omega)(f(m)<g(m)).\]

\noindent $3\Rightarrow 1$\ \ \ Given $k>0$. Let
$X\in\iso\cap\V'$. Define
\[f(n)=\min\{m>n:|X\cap[n,m)|>2k\}\mbox{\ \ for $n\in\omega$.}\]
Since $\V\cap\baire$ is unbounded in $\V'\cap\baire$ we find an increasing
function $g\in\V\cap\baire$ such that $(\exists^\infty
n\!\in\!\omega)(f(n)<g(n))$. Let $k_n\in\omega$ be defined by: 
\[k_0=0,\quad k_{n+1}=k+1+k_n+g(k_n).\]
Put $K_n=[k_n,k_{n+1})$. Clearly $\langle
K_n:n\in\omega\rangle\in\pk\cap\V$. Now suppose that $m\in K_n$ is such that
$f(m)<g(m)$. As $g$ is increasing we have $g(m)<g(k_{n+1})<k_{n+2}$.
Consequently $|[m,k_{n+2})\cap X|>2k$ and hence either $|K_n\cap X|>k$ or
$|K_{n+1}\cap X|>k$. Hence $(\exists^\infty n\!\in\!\omega)(|K_n\cap X|>k)$.
\QED

In a similar way one can prove the analogous result for the
$c\rexists$-localization.
\begin{proposition}
\label{dominujaca}
\hspace{0.1in}Let $\V\subseteq\V'$ be models of ZFC. Then the following
conditions are equivalent: 
\begin{enumerate}
\item For each $k\in\omega$ the pair $(\V,\V')$ has the
$c\rexists$-localization. 
\item For some $k\in\omega$ the pair $(\V,\V')$ has the
$c\rexists$-localization.
\item $\V\cap\baire$ is a dominating family in $\V'\cap\baire$.\QED
\end{enumerate}
\end{proposition}
For the $\rexists$-localization we did not find a full description. First
note that the requirement that members of the partition have to have at
least $k+2$ elements (i.e.  $\rexists\subseteq\icso\times\P_{k+1}$) is to
avoid a trivial localization. If we divide $\omega$ into $k+1$-element
intervals then for each set $X\in\icso$ infinitely many intervals contain at 
most $k$ members of $X$.

\begin{definition}
Let $\V\subseteq\V'$ be models of ZFC and let $k<l<\omega$. A set
$X\in\icso\cap\V'$ is called {\em $(l,k)$-large (over $\V$)} 
if for every sequence $\langle K_n: n\in\omega\rangle\in\V$ of disjoint
$l$-element subsets of $\omega$ we have: $(\forall^\infty n\in\omega)(|\;
K_n\cap X\;|>k)$.

\noindent [Note that we do not require that $\langle K_n: n\in\omega\rangle$
is a partition of $\omega$.]
\end{definition}

\begin{theorem}
\label{exists}
All $\rexists$-localizations (for $k\in\omega$) are equivalent, i.e. if an
extension $\V\subseteq\V'$ has the $\rexists$-localization property for some
$k\in\omega$ then it has the $\rexists$-localization for each $k\in\omega$. 
\end{theorem}

\Proof Let $\V\subseteq\V'$ be models of ZFC. Let $k<l<\omega$. 

\begin{claim}
\label{cl1}
If $X\in\icso\cap\V'$ is $(l,k)$-large,
$k+1<l<\omega$ and $m\geq 1$ then $X$ is $(lm,lm-(l-k))$-large.
\end{claim}

\noindent{\em Proof of Claim:}\ \ \ Let $\langle K_n: n\in\omega\rangle\in\V$
be a sequence of disjoint subsets of $\omega$ of the size $lm$. Let
$K_n=\{a_{n,i}: i<lm\}$ be the increasing enumeration and let
$K^A_n=\{a_{n,i}:i\in A\}$ for $A\subseteq lm$. Fix a set $A\in
[lm]^{\textstyle l}$ and consider the sequence $\langle K_n^A: n\in\omega
\rangle\subseteq [\omega]^{\textstyle l}$. This sequence is 
in $\V$ and its members are disjoint. Thus we find $N(A)$ such that
$(\forall n\geq N(A))(|\;K^A_n\cap X\;|>k)$. Let $N=\max\{N(A):A\in
[lm]^{\textstyle l}\}$. Then for each $n\geq N$ and $A\in 
[lm]^{\textstyle l}$ we have $|\;K^A_n\cap X\;|>k$ and hence
$|\;K^A_n\setminus X\;|<l-k$. Hence we conclude that $|\;K_n\setminus
X\;|<l-k$ for each $n\geq N$ (just take a suitable $A$). Thus $|\;K_n\cap
X\;|>lm-(l-k)$ for each $n\geq N$ and the claim is proved.

\begin{claim}
\label{cl2}
If $X\in\icso\cap\V'$ is $(l,k)$-large,
$k+1<l<\omega$ and $m\geq l-k$ then $X$ is $(m, m-(l-k))$-large.
\end{claim}

\noindent{\em Proof of Claim:}\ \ \ Let $\langle K_n: n\in\omega\rangle\in\V$
be a sequence of disjoint $m$-element sets. Put $K_n^*=K_{ln}\cup
K_{ln+1}\cup\ldots\cup K_{ln+(l-1)}$. Clearly $\langle K_n^*:
n\in\omega\rangle\in\V$, $|K^*_n|=lm$ and the sets $K^*_n$'s are disjoint. It
follows from \ref{cl1} that the set $X$ is $(lm,lm-(l-k))$-large and hence
there is $N$ such that 
\[(\forall n\geq N)(|\;K^*_n\cap X\;|>lm-(l-k)).\]
So $|\;K^*_n\setminus X\;|<l-k$ for $n\geq N$ and hence $|\;K_n\setminus
X\;|<l-k$ for $n\geq lN$. This implies $(\forall n\geq lN)(|\;K_n\cap
X\;|>m-(l-k))$ and the claim is proved. 

\begin{claim}
\label{cl3}
Assume that the extension $\V\subseteq\V'$
has the $R^\exists_0$-localization property. Then it has the
$\rexists$-localization for each $k\in\omega$.
\end{claim}

\noindent{\em Proof of Claim:}\ \ \ Let $k>0$ and assume that the
$\rexists$-localization fails. Then we have a set $X\in\icso$ witnessing it,
i.e. such that 
\[(\forall\langle K_n:n\in\omega\rangle\in\P_{k+1}\cap\V)(\forall^\infty
n\in\omega)(|\;K_n\cap X\;|>k).\]
Then, in particular, the set $X$ is $(l,k)$-large for each $l>k+1$. By claim
\ref{cl2} it is $(l,0)$-large for each $l\geq 2$. By the 
$R^\exists_0$-localization we find a partition $\langle K_n:
n\in\omega\rangle\in\P_1\cap\V$ such that $(\exists^\infty
n\in\omega)(K_n\cap X=\emptyset)$. Each set $K_n$ we partition into 2- and
3-element sets:
\[K_n=\bigcup\{K^2_{n,i}: i\in w^2_n\}\cup\bigcup\{K^3_{n,i}:
i\in w^3_n\}\]
(everything should be done in $\V$, of course). Next look at
\[\langle K^2_{n,i}: n\in\omega, i\in w^2_n\rangle,\quad\quad \langle
K^3_{n,i}:n\in\omega, i\in w^3_n\rangle\in\V.\]
These are sequences of disjoint 2- (or 3-, respectively) element subsets of
$\omega$. At least one of them is infinite; for simplicity we assume that both
are infinite. As $X$ is both $(2,0)$- and $(3,0)$-large we find $N$ such that
for each $n\geq N$, $j\in\{2,3\}$ and $i\in w^j_n$ we have $K^j_{n,i}\cap
X\neq\emptyset$. But this implies that for each $n\geq N$ the intersection
$K_n\cap X$ is not empty - a contradiction to the choice of $\langle K_n:
n\in\omega\rangle$. 

\begin{claim}
\label{cl4}
Assume that $k\in\omega$, $\V\subseteq\V'$ has
the $\rexists$-localization property. Then the extension $\V\subseteq\V'$
has the $R^\exists_0$-localization.
\end{claim}

\noindent{\em Proof of Claim:}\ \ \ Suppose that the $R^\exists_0$-localization
fails and this is witnessed by a set $X\in\icso\cap\V'$. As earlier we
conclude from this that the set $X$ is $(l,0)$-large and hence, by \ref{cl2},
it is $(l+k,k)$-large for each $l\geq 2$. By the $\rexists$-localization we
find $\langle K_n: n\in\omega\rangle\in\P_{k+1}\cap\V$ such that
$(\exists^\infty n\in\omega)(|\;K_n\cap X\;|\leq k)$. For $i\in\omega$ let
$A_i=\{n\in\omega: |K_n|=i\}$ (some of these sets can be finite or even
empty). For $n\notin\bigcup_{i\leq 2k+2} A_i$ partition $K_n$ into 2- and
3-element sets to have more than $k$ pieces: 
\[K_n=\bigcup\{K^2_{n,i}: i\in w^2_n\}\cup\bigcup\{K^3_{n,i}: i\in
w^3_n\},\quad\quad |w^2_n|+|w^3_n|>k\]
(everything is done in $\V$, of course). Consider the sequences
\begin{quotation}
$\langle K^2_{n,i}:i\in w^2_n, n\notin\bigcup_{j\leq 2k+2}A_j\rangle$, 

$\langle K^3_{n,i}:i\in w^3_n, n\notin\bigcup_{j\leq 2k+2}A_j\rangle$, 

$\langle K_n: n\in A_{k+2}\rangle$,\ldots,$\langle K_n: n\in A_{2k+2}\rangle$
\end{quotation}
(note that 
$A_i=\emptyset$ for $i<k+2$). These are sequences of disjoint sets of the
sizes $2,3,k+2,\ldots,2k+2$, respectively, and all sequences are in $\V$.
Since $X$ is $(l+k,k)$-large for each $l\geq 2$ and it is (2,0)- and
(3,0)-large we find $N$ such that
\begin{description}
\item[(a)] if $n>N$, $n\in A_i$, $k+2\leq i\leq 2k+2$ then $|K_n\cap X|>k$
\item[(b)] if $n>N$, $n\notin\bigcup_{j\leq 2k+2}A_j$, $i\in w^x_n$,
$x\in\{2,3\}$ then $K^x_{n,i}\cap X\neq\emptyset$.
\end{description}
The condition (b) implies that if $n>N$, $n\notin\bigcup_{j\leq 2k+2}A_j$
then $|\;K_n\cap X\;|>k$ (recall that we have more than $k$ sets
$K^x_{n,i}$). Consequently $|\;K_n\cap X\;|>k$ for all $n>N$ which
contradicts the choice of $\langle K_n: n\in\omega\rangle$.\QED

The above proof suggests to consider $(m,0)$-large sets (over $\V$) and ask
if the existence of such sets depends on $m\geq 2$. The answer is given by
the next result.

\begin{proposition}
Suppose $\V\subseteq\V'$ are models of ZFC, $m\geq 2$. Then there exists an
$(m,0)$-large set over $\V$ if and only if there exists an $(m+1,0)$-large
set over $\V$.
\end{proposition}

\Proof Clearly each $(m,0)$-large set is $(m+1,0)$-large. So suppose now that
$X\in\icso\cap\V'$ is $(m+1,0)$-large over $\V$. If it is $(m,0)$-large then
we are done. So assume that $X$ is not $(m,0)$-large and this is witnessed by
$\langle K_n:n\in\omega\rangle\in\V$ (so $|K_n|=m$, $K_n$'s are disjoint and
$(\exists^\infty n\in\omega)(K_n\cap X=\emptyset)$). Let 
\[Y=\{n\in\omega:K_n\cap X\neq\emptyset\}.\] 
Clearly $Y$ is infinite co-infinite. We are going to show that $Y$ is
$(2,0)$-large (and hence $(m,0)$-large) over $\V$. 

Suppose that $\langle L_n: n\in\omega\rangle\in\V$ is a sequence of disjoint
2-element sets. Let $K^*_n=\bigcup_{l\in L_n}K_l$. Thus $|K^*_n|=2m$ and
$K^*_n$'s are disjoint. Obviously the sequence $\langle
K^*_n:n\in\omega\rangle$ is in $\V$. Since the set $X$ is $(m+1,0)$-large (and
hence, by \ref{cl2}, $(2m,m-1)$-large) we have 
\[(\forall^\infty n\in\omega)(|\;K_n^*\cap X\;|>m-1).\]
But $K^*_n\cap X\neq\emptyset$, $L_n=\{l^0,l^1\}$ imply that either
$K_{l^0}\cap X\neq\emptyset$ or $K_{l^1}\cap X\neq\emptyset$ and hence
$L_n\cap X\neq\emptyset$. Consequently $Y$ is (2,0)-large. \QED

\begin{corollary}
Let $\V\subseteq\V'$ be models of ZFC, $m\geq 2$, $k\in\omega$. Then the
following conditions are equivalent:
\begin{enumerate}
\item there is no $(m,0)$-large set in $\icso\cap\V'$ over $\V$
\item there is no $(2,0)$-large set in $\icso\cap\V'$ over $\V$
\item $\V\subseteq\V'$ has the $R^\exists_0$-localization property
\item $\V\subseteq\V'$ has the $\rexists$-localization property.\QED
\end{enumerate}
\end{corollary}

\noindent{\sc Remark:}\ \ \ One can consider a modification of the notion of
$(m,0)$-largeness giving (probably) more freedom. For an increasing function
$f\in\baire\cap\V$ we say that a set $X\in\icso$ is $f$-large over $\V$ if
for every sequence $\langle K_n: n\in\omega\rangle\in\V$ of disjoint finite
subsets of $\omega$ we have
\[\mbox{either }(\exists n\in\omega)(|K_n|<f(n)+2)\quad\mbox{ or }
(\exists^\infty n\in\omega)(|K_n\cap X|>f(n)).\]

\begin{proposition}
\label{meabou}
Let $\V\subseteq\V'$ be models of ZFC.\\ 
{\bf a)}\ \ \ If $\V\cap\can$ is not meager iv $\V'$ then the pair
$(\V,\V')$ has the $R^\exists_0$-localization property.\\ 
{\bf b)}\ \ \ If the pair $(\V,\V')$ has the $R^\exists_0$-localization
property then $\V\cap\baire$ is unbounded in $\V'\cap\baire$.
\end{proposition}

\Proof {\bf a)}\ \ If $\V\cap\can$ is not meager in $\V'$ then 
$$(\forall f\!\in\!\V'\!\cap\!\baire)(\forall
Y\!\in\!\V'\!\cap\!\iso)(\exists g\!\in\!\V\!\cap\!\baire)(\exists^\infty
n\!\in\!Y)(f(n)=g(n))\leqno(*)$$  
and $\V\cap\baire$ is unbounded in $\V'\cap\baire$ (see \cite{Bar1}). By
proposition \ref{unbounded}, the pair $(\V,\V')$ has the
$cR^\forall_1$-localization property. 
Suppose that $X\in\icso\cap\V'$. By the $cR^\forall_1$-localization we find
a partition $\langle K_n:n\in\omega\rangle\in\V\cap\P_1$ such that
$(\exists^\infty n\in\omega)(|\;K_n\setminus X\;|\geq 2)$. In $\V'$ we
define 
\begin{quotation}
$f(n)=K_n\setminus X\in\fsuo$\ \ \ \ (for $n\in\omega$),

$Y=\{n\in\omega:|f(n)|\geq 2\}$. 
\end{quotation}
By $(*)$ we find $g\in\V$, $g:\omega\lra\fsuo$ such that $g(n)\in
[K_n]^{\textstyle\geq 2}$ and 
\[(\exists^\infty n\in Y)(f(n)=g(n)).\]
Since $f(n)=g(n)$ implies $g(n)=K_n\setminus X$ we get
$(\exists^\infty n\in\omega)(g(n)\cap X=\emptyset)$. Hence we easily
get that $X$ can be $R^\exists_0$-localized by a partition from $\V$.

{\bf b)} Since $|K|>1\ \&\ K\cap(\omega\setminus X)=\emptyset$ implies
$|K\cap X|\geq 2$ we get that $R^\exists_0$-localization implies the
$cR^{\forall}_1$-localization. Now proposition~\ref{unbounded} works.\QED

The next result gives some bounds on possible improvements of the previous
one. 

\begin{proposition}
\label{cora}
\begin{enumerate}
\item The Cohen forcing notion has the $R^\exists_0$-localization property.
Consequently, the $R^\exists_0$-localization does not imply that the old
reals are a dominating family.
\item The Random real forcing does not have the $R^\exists_0$-localization
property. Consequently, the localization is not implied by the fact that
there is no unbounded real in the extension.
\end{enumerate}
\end{proposition}

\Proof {\em 1.}\ \ \ As in the extensions via the Cohen forcing the
ground model reals are not meager, we may apply \ref{meabou}.

{\em 2.}\ \ \ The Random algebra $\random$ is the quotient algebra of Borel
subsets of $\can$ modulo the ideal of Lebesgue null sets. We define a
$\random$-name for an element of $\icso$ that cannot be localized:
\begin{quotation}
\noindent Let $l_0=0$, $l_{k+1}=l_k+2^{k^2}$ (for $k\in\omega$). 

\noindent For each $k\in \omega$ fix disjoint Borel sets
$A_m\subseteq\can$ for $l_k\leq m<l_{k+1}$ such that $\mu(A_m)=2^{-k^2}$,
where $\mu$ is the Lebesgue measure on $\can$.  

\noindent $\dot{X}$ is a $\random$-name for a subset of $\omega$ such that 
if $m\in [l_k,l_{k+1})$, $k\in\omega$ then $\lbv m\notin\dot{X}\rbv_{\random}=
[A_m]_{\mu}$.   
\end{quotation}
It should be clear that $\forces_{\random}(\forall
k\in\omega)(|\;[l_k,l_{k+1})\setminus\dot{X}\;|=1)$. 

\noindent Suppose now that $\langle K_n:n\in\omega\rangle\in{\cal P}_1\cap\V$.
Let $k^0_n=\min K_n$ and $k^1_n=\max K_n$. Suppose that $l_m\leq
k^0_n<l_{m+1}$. Note that 
\[k^1_n<l_{m+1}\ \ \Rightarrow\ \ \lbv k^0_n\notin \dot{X}\
\&\ k^1_n\notin\dot{X}\rbv_{\random}={\bf 0}\]
and 
\[k^1_n\geq l_{m+1}\ \ \Rightarrow\ \ \mu(\lbv
k^1_n\notin\dot{X}\rbv_{\random})\leq 2^{-(m+1)^2}.\]
Hence 
\[\mu(\lbv\dot{X}\cap K_n=\emptyset\rbv_{\random})\leq\mu(\lbv
k^0_n\notin\dot{X}\ \&\ k^1_n\notin\dot{X}\rbv_{\random})\leq
2^{-(m+1)^2}.\] 
Consequently, for each $m\geq 0$:
\[\mu(\lbv(\exists n\in\omega)(l_m\leq k^0_n\ \&\
K_n\cap\dot{X}=\emptyset)\rbv_{\random})\leq\sum_{n\in\omega}\mu(\lbv l_m\leq
k^0_n\ \&\ K_n\cap\dot{X}=\emptyset\rbv_{\random})=\]
\[=\sum_{r\geq m}\big(\sum_{k^0_n\in [l_r,l_{r+1})} \mu(\lbv\dot{X}\cap
K_n=\emptyset\rbv_{\random})\big)\leq\sum_{r\geq m}
2^{r^2}2^{-(r+1)^2}=\frac{1}{3} 2^{1-2m}.\]
Hence we can conclude that $\mu(\lbv(\exists^\infty
n\in\omega)(K_n\cap\dot{X} =\emptyset)\rbv_{\random})=0$ which means
$\forces_{\random}$``$\langle K_n: n\in\omega\rangle$ does not
$R^\exists_0$-localize $\dot{X}$''.\QED
\medskip

Though the random real forcing is an example of a forcing notion adding a
(2,0)-large set over $\V$ (without adding an unbounded real!) it does not
seem to be the minimal one. A canonical example of a forcing notion without
the $R^\exists_0$-localization property is given below. (Recall that a forcing
notion $\q$ is $\sigma$-centered if it can be presented as a countable union of
sets which all finite subsets have upper bounds in $\q$.)

\begin{example}
There is a $\sigma$-centered (Borel) forcing notion $\q$ adding no
dominating real and without the $R^\exists_0$-localization property.
\end{example}

\Proof The forcing notion $\q$ consists of pairs $(u,\K)$ such that
$u\in\fsuo$ and $\K$ is a finite set of families of disjoint 2-element
subsets of $\omega$ (so $F\in\K\ \Rightarrow\ F\subseteq[\omega]^{\textstyle
2}$). The order of $\q$ is given by
\begin{quotation}
\noindent $(u_0,\K_0)\leq (u_1, \K_1)$ if and only if 

\noindent $u_1\cap (1+\max u_0)=u_0$, $\K_0\subseteq\K_1$ and if $K\in F\in
\K_0$ and $K\subseteq u_1$ then $K\subseteq u_0$.
\end{quotation}
For $u\in\fsuo$, $m\in\omega$ let 
\[Q_u=\{(u,\K): (u,\K)\in\q\}\ \mbox{ and }\ Q^m_u=\{(u,\K)\in Q_u:
|\K|=m\}.\]
Since each $Q_u$ is obviously centered we get that $\q$ is $\sigma$-centered.
Let $\dot{w}$ be a $\q$-name such that $\dot{w}^G=\bigcup\{u:
(\exists\K)((u,\K)\in G)\}$ for each generic filter $G\subseteq\q$ over $\V$.
Note that $\forces_{\q}\dot{w}\in\icso$. If $K\in F\in\K$, $(u,\K)\in\q$ and
$\max u<\min K$ then $(u,\K)\forces_{\q} K\setminus\dot{w}\neq\emptyset$.
Hence we conclude that 
\[\forces_{\q}\mbox{``}\omega\setminus\dot{w}\mbox{ is } (2,0)\mbox{-large
over }\V\mbox{''}.\]
So $\q$ does not have the $R^\exists_0$-localization property. Suppose now
that $\tau$ is a $\q$-name for a member of $\baire$.
\begin{claim}
\label{cl5}
Let $u\in\fsuo$, $m,n\in\omega$. The there is
$f(u,m,n)<\omega$ such that:
\begin{quote}
\noindent for every $p\in Q^m_u$ there is $q=(u^q,\K^q)\geq p$ such that $q$
decides the value of $\tau(n)$ and $\max u^q<f(u,m,n)$.
\end{quote}
\end{claim}

\noindent{\em Proof of Claim:}\ \ \ The space $\X$ of all families
$F\subseteq [\omega]^{\textstyle 2}$ of two-element disjoint subsets of
$\omega$ can be equipped with a natural topology. For $F\in\X$, $N\in\omega$
the $N$-th basic open neighbourhood of $F$ is 
\[\{F'\in\X: \{K\cap N: K\in F\}=\{K\cap N: K\in F'\}\}.\] 
This topology is compact. It introduces a (product) topology on $Q^m_u$ such
that if $q\in Q_{u'}$, $q\geq p$, $p\in Q^m_u$ then for some open
neighbourhood $V$ of $p$ (in $Q^m_u$) each member of $V$ has an extension in
$Q_{u'}$. Applying this fact and the compactness of $Q^m_u$ we get the claim. 

Now we define a function $g\in\baire$ putting 
\[\begin{array}{ll}
g(k)=1+\max\{l: &(\exists u\subseteq k)(\exists m,n\leq k)(\exists
v)\\
\ &(u\subseteq v\subseteq f(u,m,n)\ \mbox{ and }\ (\exists q\in Q_v)(q\forces
\tau(n)=l))\}.\\
\end{array}\] 
Given $p\in Q^m_u$, $l\in\omega$. Take $k>\max\{l,m,\max u\}$. By the
definition of $f(u,m,k)$ we find $v$ such that $u\subseteq v\subseteq
f(u,m,k)$ and some $q\in Q_v$, $q\geq p$ decides the value of $\tau(k)$. By
the definition of the function $g$, the condition $q$ forces
``$\tau(k)<g(k)$''.\QED 
\medskip

\noindent{\sc Remark:}\ \ \ An example of a forcing notion $\p$ with the
$R^\exists_0$-localization property and such that 
\[\forces_{\p}\mbox{``}\V\cap\can\mbox{ is meager''}\]
is an application of a general framework of \cite{RoSh:470} and will be
presented there. 

\section{Between dominating and unbounded reals}

In this section we are interested in some localizations which are between the
$c\rall$-localization and $c\rexists$-localization (so between not adding a
dominating real and not adding an unbounded real).  The localizations are
similar to that considered in the previous section. The difference is that
we will consider partitions of $\omega$ into intervals and we will introduce
quantifiers stronger than $\exists^\infty n$ but weaker than $\forall^\infty
n$.

For an infinite subset $X$ of $\omega$ let $\mu_X:\omega\lra X$ be the
increasing enumeration of $X$. A set $X\in\iso$ can be identified with the
partition 
\[\langle [\mu_X(n),\mu_X(n+1)): n\in\omega\rangle\]
of $\omega\setminus\mu_X(0)$, so essentially $\iso\subseteq\P_0$. 
Now, for $k>0$ and an increasing function $\phi\in\baire$, we define relations
$S_k, S_+, S_{+\epsilon}, S_+^\phi\subseteq\iso\times\iso$: 
\[(X,Y)\in S_k\equiv (\exists^\infty n\!\in\!\omega)(\forall
i<k)(|\;[\mu_Y(n+i), \mu_Y(n+i+1))\cap X\;|\geq 2)\]
\[(X,Y)\in S_+\equiv (\forall m\!\in\!\omega)(\exists n\!\in\!\omega)(\forall
i<m)(|\;[\mu_Y(n+i), \mu_Y(n+i+1))\cap X\;|\geq 2)\]
\[(X,Y)\in S_{+\epsilon}\equiv (\exists^\infty n\!\in\!\omega)(\forall
i<2^n)(|\;[\mu_Y(2^n+i), \mu_Y(2^n+i+1))\cap X\;|\geq 2)\]
\[(X,Y)\in S_+^\phi\equiv (\exists^\infty n\!\in\!\omega)(\forall
i<\phi(n))(|\;[\mu_Y(n+i), \mu_Y(n+i+1))\cap X\;|\geq 2).\] 
[The relation $S_{+\epsilon}$ appears here for historical reasons only: it
determines a cardinal invariant which was of serious use in \cite{RoSh:475}.]

\noindent Note that $S_1\supseteq S_2\supseteq S_3\supseteq\ldots\supseteq
S_+\supseteq S_{+\epsilon}\cup S_+^\phi$ (remember that $\phi$ is increasing).
If the function $\phi$ is increasing fast enough (e.g. $\phi(n)>2^{2n}$) then
$S_{+\epsilon}\supseteq S_+^\phi$. 
\smallskip 

\noindent It should be clear that if we consider $S_*$-localizations we could
put any integer greater than 2 in place of 2 in the definitions above. However
we do not know if replacing 2 by 1 provides the same notions of localizations. 

\begin{proposition}
\label{ogolne}
\hspace{0.1in} Let $\V\subseteq\V'$ be models of ZFC.
\begin{description}
\item[(a)] The pair $(\V,\V')$ has the $S_1$-localization
property if and only if $\V\cap\baire$ is unbounded in $\V'\cap\baire$.
\item[(b)] If the pair $(\V,\V')$ has the
$S_{k+1}$-localization property then it has the $S_k$-localization. The
$S_{+}$-localization property implies the $S_k$-localization for each $k>0$
and is implied by both the $S_+^\phi$ and the $S_{+\epsilon}$-localization
properties ($\phi$ - an increasing function).
\item[(c)] If $\V\cap\baire$ is dominating in $\V'\cap\baire$ then the
extension $\V\subseteq\V'$ has the $S_+^\phi$-localization property for
every increasing function $\phi\in\V\cap\baire$.\QED
\end{description}
\end{proposition}

\begin{proposition}
\label{non}
If $\V\subseteq\V'$ are models of ZFC such that $\V\cap\can$ is not meager in
$\V'$, $\phi\in\V\cap\baire$ is increasing then the pair $(\V,\V')$ has the
$S_+^\phi$-localization property.
\end{proposition}

\Proof The proof is almost the same as that of \ref{meabou}(a). 
Suppose that $X\in\V'\cap\iso$. Let $X_0\in [X]^{\textstyle \omega}$ be such
that for each $n\in\omega$
\[|\;[\mu_{X_0}(n),\mu_{X_0}(n+1))\cap X\;|>3\phi(\mu_{X_0}(n)+4)+6.\]
As $\baire\cap\V$ is unbounded in $\baire\cap\V'$ we find a set
$Y_0\in\V\cap\baire$ such that the set
\[X_1=\{\mu_{Y_0}(n): n\in\omega\ \&\ 2\leq
|\;[\mu_{Y_0}(n),\mu_{Y_0}(n+1))\cap X_0\;|\}\]
is infinite. Let $f:\omega\lra\fsuo$ be such that for every $\mu_{Y_0}(n)\in
X_1$ we have
\[f(\mu_{Y_0}(n))=[\mu_{Y_0}(n),\mu_{Y_0}(n+1))\cap X\]
(so, for $\mu_{Y_0}(n)\in X_1$, $|\;f(\mu_{Y_0}(n))\;|>
3\phi(\mu_{X_0}(k)+4)+6$ where $k$ is the first such that $\mu_{Y_0}(n)\leq
\mu_{X_0}(k)<\mu_{Y_0}(n+1)$). By $(*)$ from the proof of \ref{meabou}
we find a function $g\in\V$, $g:\omega\lra\fsuo$ such that
\[(\exists^\infty n\in\omega)(\mu_{Y_0}(n)\in X_1\ \ \&\ \
g(\mu_{Y_0}(n))=f(\mu_{Y_0}(n))).\] 
Next, using $g$ and $Y_0$ (both are in $\V$) we define sets
$Y_1,Y\in\V\cap\iso$:
\[Y_1=\bigcup_{n\in\omega}g(\mu_{Y_0}(n))\cap [\mu_{Y_0}(n),\mu_{Y_0}(n+1))\]
and $\mu_Y(n)=\mu_{Y_1}(3n)$ for each
$n$.  Note that if $n$ is such that $\mu_{Y_0}(n)\in X_1$ and
$g(\mu_{Y_0}(n))=f(\mu_{Y_0}(n))$ then
\[Y_1\cap [\mu_{Y_0}(n),\mu_{Y_0}(n+1)) = X\cap [\mu_{Y_0}(n),\mu_{Y_0}(n+1))\]
is of the size $> 3\phi(\mu_{Y_0}(n)+4)+6$ and consequently $Y$
$S_{+}^\phi$-localizes $X$. \QED   

\begin{proposition}
\label{cos}
\begin{enumerate}
\item The Cohen forcing notion has the $S^\phi_+$-localization property for
each increasing function $\phi\in\baire$. 
\item If $\V\subseteq\V'\subseteq\V''$ are models of ZFC, the pair $(\V,\V')$
has the $S^\phi_+$-localization property and $\V'\cap\baire$ is dominating in
$\V''\cap\baire$ then the extension $\V\subseteq\V''$ has the
$S^\phi_+$-localization property.  
\item The iteration of the Cohen forcing notion and the random real forcing
has the $S^\phi_+$-localization.
\end{enumerate}
\noindent Consequently, the $S^\phi_+$-localization property implies neither
that the ground model reals are dominating in the extension nor that the old
reals are not meager.\QED
\end{proposition}

\begin{theorem}
\label{eska}
For each $k>0$, the $S_k$-localization property does not imply the
$S_{k+1}$-localization property. 
\end{theorem}

\Proof To prove the theorem we will define a forcing notion $\q_k$ which
will possess the $S_k$-localization property but not $S_{k+1}$. The forcing
notion is similar to that used in \cite{Sh:207}, \cite{BsSh:242} and is a
special case of the forcing notions of \cite{RoSh:470}.  

We start with a series of definitions.

\begin{enumerate}
\item A function ${\bf n}$ is a {\em nice norm on $A$} if:
\begin{itemize}
\item ${\bf n}:\P(A)\lra\omega$ and it is monotonic 

(i.e. $B\subseteq C\subseteq A\ \ \Rightarrow\ \ {\bf n}(B)\leq {\bf n}(C)$),
\item if $B\subseteq C\subseteq A$, ${\bf n}(C)>0$ 

then either ${\bf n}(B)\geq{\bf n}(C)-1$ or ${\bf n}(C\setminus B)\geq{\bf
n}(C)-1$,  
\item ${\bf n}(A)>0$, if $a\in A$ then ${\bf n}(\{a\})\leq 1$.
\end{itemize}
\item {\em A creature} is a tuple $\t=\langle T,\nor,L,R\rangle$ such that 
\begin{description}
\item[$\alpha)$] $T\subseteq\fseo$ is a finite nonempty tree,
\item[$\beta)$] for each $t\in T$, either $\suc_T(t)=\emptyset$ or
$|\suc_T(t)|=k$ or $|\suc_T(t)|> k$

[so we have three kinds of nodes in the tree $T$]
\item[$\gamma)$] $\nor$ is a function with the domain 
\[\dom(\nor)=\{t\in T:|\suc_T(t)|>k\}\]
and such that for $t\in\dom(\nor)$, $\nor(t)$ is a nice norm on $\suc_T(t)$

[$\nor$ stands for ``norm''],
\item[$\delta)$] if $s\in\suc_T(t)$ then either $|\suc_T(t)|\neq k$ or
$|\suc_T(s)|\neq k$

[i.e. we do not have two successive $k$-ramifications in $T$],
\item[$\varepsilon)$] $L, R: T\lra\omega$ are functions such that for each
$t\in T$: 
\begin{itemize}
\item $L(t)\leq R(t)$
\item if $s\in\suc_T(t)$ then $[L(s),R(s)]\subseteq[L(t),R(t)]$
\item if $s_1,s_2\in\suc_T(t)$ are distinct\\
then $[L(s_1),R(s_1)]\cap [L(s_2), R(s_2)]=\emptyset$
\item if $\suc_T(t)=\emptyset$ then $L(t)=R(t)$
\end{itemize}
[$L$ stands for ``left'' and $R$ is for ``right''].
\end{description}
\item We will use the convention that if $a,b$ are indexes and $\t^a_b$ is a
creature then its components are denoted by $T^a_b$, $\nor^a_b$, $L^a_b$ and 
$R^a_b$, respectively.
\item Let $\t=\langle T,\nor,L,R\rangle$ be a creature. We define its weight
$\|\t\|$ and its contribution $\cont(\t)$:
\[\|\t\|=\min\{\nor(t)(\suc_T(t)): t\in\dom(\nor)\ \&\ |\suc_T(t)|>k\}\]
[if $\dom(\nor)=\emptyset$ then we put $\|\t\|=0$]
\[\cont(\t)=\{L(t): t\in T\ \&\ \suc_T(t)=\emptyset\}\]
[recall that if $t$ is a leaf in $T$ then $L(t)=R(t)$].
\item Let $\t_i$ ($i=0,1$) be creatures. We say that the creature $\t_1$ {\em
refines} $\t_0$ (we write: $\t_0\leq\t_1$) if:
\begin{description}
\item[{\bf a)}] $T_1\subseteq T_0$, $L_1=L_0\rest T_1$, $R_1=R_0\rest T_1$,
\item[{\bf b)}] if $t\in T_1$ then 
\begin{quotation}
$|\suc_{T_1}(t)|>k$ iff $|\suc_{T_0}(t)|>k$, 

$|\suc_{T_1}(t)|=k$ iff $|\suc_{T_0}(t)|=k$, and

$|\suc_{T_1}(t)|=\emptyset$ iff $|\suc_{T_0}(t)|=\emptyset$
\end{quotation}
(in other words we keep the kind of nodes),
\item[{\bf c)}] $\nor_1(t)=\nor_0(t)\rest\P(\suc_{T_1}(t))$ for all
$t\in\dom(\nor_1)$. 
\end{description}
\item Let $\t_0,\t_1,\ldots,\t_n$ be creatures. We say that the creature
$\t$ is built of the creatures $\t_0,\ldots,\t_n$
(we will write it as $\t\in\Sigma(\t_0,\ldots,\t_n)$) if there is a maximal
antichain $F$ of $T$ such that for each $t\in F$, for some $i\leq n$:
\begin{itemize}
\item $\{s\in T:t\subseteq s\}=\{t\hat{\ }r: r\in T_i\}$,
\item if $s=t\hat{\ }r\in T$, $r\in T_i$ then $L(s)=L_i(r)$,
$R(s)=R_i(r)$ and $\nor(s)=\nor_i(r)$ (if defined). 
\end{itemize}
\item If $n\geq k$, $H:\P(n+1)\lra\omega$ is a nice norm on $n+1$ (i.e. it
satisfies the conditions listed in (1)) then $S_H(\t_0,\ldots,\t_n)$ is the
creature $\t$ such that 
\[T=\{\langle i\rangle\hat{\ }t: i\leq n, t\in T_i\},\quad
\nor(\langle\rangle)=H,\quad \nor(\langle i\rangle\hat{\ }t)=\nor_i(t)\]
and similarly for $L$ and $R$.

Clearly $S_H(\t_0,\ldots,\t_n)\in \Sigma(\t_0,\ldots,\t_n)$. 
\item We define gluing $k$ creatures similarly to the operation $S_H$ above. 
Thus $S(\t_0,\ldots,\t_{k-1})=\t$ is a creature such that 
\[T=\{\langle i\rangle\hat{\ }t: i<k, t\in T_i\}\]
and $\nor$, $L$, $R$ are defined naturally. Once again,
$S(\t_0,\ldots,\t_{k-1})\in \Sigma(\t_0,\ldots,\t_{k-1})$.  
\item For a creature $\t$ we define its upper half $\t^\uh=\langle
T^\uh,\nor^\uh,L^\uh,R^\uh\rangle$ by: 
\[T^\uh=T,\quad L^\uh=L,\quad R^\uh=R\quad\mbox{ but}\]
\[\nor^\uh(t)(A)=\max\{0,\nor(t)(A)-[\frac{\|\t\|}{2}]\}\]
whenever $t\in T$, $|\suc_T(t)|>k$ and $A\subseteq\suc_T(t)$. Above, $[x]$
stands for the integer part of $x$. 

[It is routine to check that $\t^\uh$ is really a creature and that
$\|\t^\uh\|=\|\t\|-[\frac{1}{2}\|\t\|]$, $\cont(\t^\uh)=\cont(\t)$.]
\item For creatures $\t_0,\ldots,\t_n$, the closure of $\{\t_0,\ldots,\t_n\}$
under the operations of shrinking (refining), taking the upper half and
building creatures is denoted by $\Sigma^*(\t_0,\ldots,\t_n)$. Thus 
$\{\t_0,\ldots,\t_n\}\subseteq\Sigma^*(\t_0,\ldots,\t_n)$, if
$\t\in\Sigma^*(\t_0,\ldots,\t_n)$ then 
\[\t^\uh\in\Sigma^*(\t_0,\ldots,\t_n),\quad \t\leq\t'\ \ \Rightarrow\ \
\t'\in\Sigma^*(\t_0,\ldots,\t_n),\]
and if $\t_0^\prime,\ldots,\t_m^\prime\in\Sigma^*(\t_0,\ldots,\t_n)$ then 
\[\Sigma(\t_0^\prime,\ldots,\t_m^\prime)\subseteq\Sigma^*(\t_0,\ldots,\t_n).\]
[Note that $\Sigma^*(\t_0,\ldots,\t_n)$ is finite (up to isomorphism).]
\end{enumerate}
\smallskip

Now we may define our forcing notion $\q_k$:
\medskip

\noindent{\bf Conditions}\ \ \ are sequences $\langle w,\t_0,\t_1,\t_2,\ldots
\rangle$ such that $w\in\fsuo$, $\t_i$ are creatures, 
$\|\t_i\|\lra\infty$ and 
\[\max(w)<L_0(\langle\rangle)\leq R_0(\langle\rangle)<L_1(\langle\rangle)\leq
R_1(\langle\rangle)<\ldots\]
(recall that $\t_i=\langle T_i,\nor_i, L_i, R_i\rangle$).
\medskip

\noindent{\bf The order}\ \ \ \ is given by\\
$\langle w,\t_0,\t_1\ldots\rangle\leq\langle w',\t_0',\t_1'\ldots\rangle$
if and only if\\ 
for some increasing sequence $n_0<n_1<n_2<\ldots<\omega$ 
\[w\subseteq w'\subseteq w\cup\bigcup_{i<n_0}\cont(\t_i)\quad\mbox{ and }\quad
(\forall i\in\omega)(\t_i'\in\Sigma^*(\t_{n_i},\ldots, \t_{n_{i+1}-1})).\] 
\medskip

\noindent We say that a condition $\langle w',\t_0',\t_1'\ldots\rangle\in\q_k$
is {\em a pure extension } of a condition $\langle
w,\t_0,\t_1\ldots\rangle\in\q_k$ if  
\[\langle w,\t_0,\t_1\ldots\rangle\leq\langle w',\t_0',\t_1'\ldots\rangle\quad
\mbox{ and }\quad w=w'.\]  

One can easily check that $(\q_k,{\leq})$ is a partial order and that the
relation of pure extension is transitive. The proof that the forcing notion
$\q_k$ has the required properties is broken into several claims. The first
claim is of a technical character, but it implies in particular that $\q_k$ is
proper. (The proof of this claim is straightforward and we will omit it.)

\begin{claim}
\label{cl6}
If $p\in\q_k$, $\tau_n$ (for $n<\omega$) are $\q_k$-names for ordinals then
there is a pure extension $q=\langle w,\t_0,\t_1,\ldots\rangle$ of $p$ such
that  
\begin{enumerate}
\item $\langle\|\t_i\|:i\in\omega\rangle$ is increasing,
\item for each $n$, $v\subseteq R_n(\langle\rangle)$ and $i\leq n$, 

{\em if} some pure extension of $\langle w\cup
v,\t_{n+1},\t_{n+2}\ldots\rangle$ decides the value of $\tau_i$ 

{\em then}  $\langle w\cup v,\t_{n+1},\t_{n+2}\ldots\rangle$ decides it.   
\end{enumerate}
\end{claim}

\begin{claim}
\label{cl7}
Suppose that $\t=\langle T,\nor,L,R\rangle$ is a creature, $\|\t\|\geq 15$
and $B\in\iso$. Then there is a creature $\t'\geq\t$ such that
$\|\t'\|\geq\|\t\|-14$ and  
\[(\forall n\in\omega)(\exists i\leq k)(|\;\cont(\t')\cap
[\mu_B(n+i),\mu_B(n+i+1))\;|<2).\]
\end{claim}

\noindent{\em Proof of Claim:}\ \ \ We prove this essentially by the induction
on $|T|$ (or the height of $T$). We show how to ``eliminate'' $\langle\rangle$
and apply the inductive hypothesis to $\t$ above $t$ (for
$t\in\suc_T(\langle\rangle)$).  

{\em Case 1:}\ \ \ $k=|\suc_T(\langle\rangle)|$

If, for each $t\in\suc_T(\langle\rangle)$, $\suc_T(t) = \emptyset$ then
$|\;\cont(\t)\;|=k$ and there is no problem. So there are $t\in
\suc_T(\langle\rangle)$ such that $|\suc_T(t)|>k$ (here we use the requirement
$(\delta)$ of the definition of creatures); above each such $t$ we can apply
the inductive hypothesis and shrink suitably the tree $T$. However the
problems coming from distinct $t$ could accumulate. Therefore for each such
$t$ we first choose a set $A=A^t\subseteq\suc_T(t)$ such that $\nor(t)(A)\geq
\|\t\|-7$ and one of the following occurs: 
\begin{itemize}
\item for some $m$, $(\forall s\in A)(\mu_B(m)\leq L(s)\leq R(s)\leq
\mu_B(m+1))$ 
\item there are $1\leq m_0<m_1$ such that 
\[(\forall s\in A)(\mu_B(m_0)\leq L(s)\leq R(s)\leq\mu_B(m_1)),\]
and if $t_0\in\suc_T(\langle\rangle)$, $t_0\neq t$ then 
\[\mbox{either }\quad R(t_0)<\mu_B(m_0-1)\quad\mbox{ or }\quad\mu_B(m_1+1)<
L(t_0).\]  
\end{itemize}
To find such a set $A$ we use the second property of nice norms. First we look
at the set
\[A_0=\{s\in\suc_T(t): (\exists m\in\omega)(\mu_B(m)\leq L(s)\leq R(s)\leq
\mu_B(m+1))\}.\] 
If $\nor(t)(A_0)\geq\nor(t)(\suc_T(t))-1\geq\|\t\|-1$ then we easily finish:
either for some $m$
\[\nor(t)(\{s\in\suc_T(t): \mu_B(m)\leq L(s)\leq R(s)\leq\mu_B(m+1)\})\geq
\|\t\|-7\]
or we have
$m_0^\prime<m_0^{\prime\prime}<m_0<m_1<m_1^\prime<m_1^{\prime\prime}$ 
such that 
\[\nor(t)(\{s\in\suc_T(t): \mu_B(m_0)\leq L(s)\leq R(s)\leq \mu_B(m_1)\})\geq
\|\t\|-7\quad\quad\mbox{ and}\]
\[(\exists s_0,s_1{\in}\suc_T(t))(\mu_B(m_0^\prime)\leq L(s_0)\leq
\mu_B(m_0^{\prime\prime})\ \&\ \mu_B(m_1^\prime)\leq R(s_1)\leq
\mu_B(m_1^{\prime\prime})).\]
So suppose that $\nor(t)(A_0)<\|\t\|-1$. In this case
$\nor(t)(\suc_T(t)\setminus A_0)\geq\|\t\|-1$. For each $s\in
\suc_T(t)\setminus A_0$ there is $m\in\omega$ such that  
\[L(s)<\mu_B(m)<R(s).\]
Removing 4 extreme points from $\suc_T(t)\setminus A_0$ (the first two and the
last two, counting according to the values of $L$) we get the required set $A$
(with $\nor(t)(A)\geq \|T\|-5$). 

\noindent Now we would like to apply the induction hypothesis above each
$t\in\suc_T(\langle\rangle)$ restricting ourselves to successors of $t$ from
the suitable set $A^t$ (if applicable). A small difficulty is that we have
decreased the norm of the creature above those $t$ (possibly by 7, as the
result of restricting to $A^t$). But now we apply the procedure described in
the case 2 below and we pass to $B^t\subseteq A^t$ such that 
\[\nor(t)(b^t)\geq\nor(t)(A^t)-7\geq\|\t\|-14.\]
Next we apply the inductive hypothesis above each $s\in B^t$,
$t\in\suc_T(\langle\rangle)\cap\dom(\nor)$. In this way we construct $\t'$ as
required (the point is that restricting to the sets $B^t\subseteq A^t$ causes
that what happens above distinct $t\in\suc_T(\langle\rangle)$ is isolated in a
sense).   
\medskip

{\em Case 2:}\ \ \ $k<|\suc_T(\langle\rangle)|$

Then one of the following possibilities occurs ($i=0,1$):
\begin{description}
\item[$(\alpha)_i$] $\nor(\langle\rangle)(\{t\in\suc_T(\langle\rangle):
(\exists m)(\mu_B(2m+i)\leq L(t)\leq R(t)<$

\hspace{3.5cm} $<\mu_B(2m+i+1))\})\geq\|\t\|-3$
\item[$(\beta)$]  $\nor(\langle\rangle)(\{t\in\suc_T(\langle\rangle):
(L(t), R(t))\cap B\neq\emptyset\})\geq\|\t\|-3$
\end{description}
If one of cases $(\alpha)_0$, $(\alpha)_1$ holds then the creature refining
$\t$ and determined by the respective set of successors of $\langle\rangle$
can serve as $\t'$. In case $(\beta)$ divide the set
\[A=\{t\in\suc_T(\langle\rangle): (L(t), R(t))\cap B\neq\emptyset\}\]
into four disjoint subsets, each containing every fourth member of $A$. One of
these subsets (call it $A^*$) has the norm $\geq\|\t\|-7$. For each $t\in A^*$
apply the inductive hypothesis to the creature given by $\{s\in T: t\subseteq
s\}$. Note that either it is of the weight $\geq\|\t\|$ or
$\suc_T(t)=\emptyset$ or $|\suc_T(t)|=k$ and for $s\in\suc_T(t)$ we have
$\suc_T(s)=\emptyset$. The last two cases are trivial and actually should be
considered separately (compare Case 1). In this way we get the required
creature $\t'$. 

\begin{claim}
\label{cl8}
$\q_k$ does not have the $S_{k+1}$-localization property.
\end{claim}

\noindent{\em Proof of Claim:}\ \ \ Let $\dot{w}$ be a $\q_k$ name such that if
$G\subseteq\q_k$ is a generic then $\dot{w}^G=\bigcup\{w: (\exists
\langle\t_0,\t_1,\ldots\rangle)(\langle w,\t_0,\t_1,\ldots\rangle\in G)\}$.
We claim that the $S_{k+1}$-localization always fails for $\dot{w}^G$, i.e.
that if $B\in\V\cap\iso$ then  
\[\forces_{\q_k}(\forall^\infty n\in\omega)(\exists i\leq k)(|\;\dot{w}\cap
[\mu_B(n+i),\mu_B(n+i+1))\;|<2).\] 
Let $p=\langle w,\t_0,\t_1,\ldots\rangle\in\q_k$ be given. We may assume that
\begin{enumerate}
\item $(\forall l\in\omega)(|\;[R_l(\langle\rangle),
L_{l+1}(\langle\rangle))\cap B\;|>2)$
\item $(\forall l\in\omega)(\|\t_l\|>15)$.
\end{enumerate}
For each $\t_i$ take the creature $\t_i^\prime\geq\t_i$ given by claim
\ref{cl7} for $\t_i$ and $B$. Look at the condition $q=\langle
w,\t^\prime_{0},\t^\prime_{1},\ldots\rangle\in \q_k$. Clearly $q\geq p$
and if $\mu_B(n)>\max(w)$ then  
\[q\forces_{\q_k}(\exists i\leq k)(|\dot{w}\cap
[\mu_B(n+i),\mu_B(n+i+1))|<2).\]
This proves the claim.
\medskip

The next claim explains why we introduced the operation of taking the upper
half of a creature as a part of the definition of (the order of) $\q_k$

\begin{claim}
\label{cl9}
Let $p=\langle w,\t_0,\t_1\ldots\rangle\in\q_k$, $m\in\omega$. Suppose that
$\tau$ is a $\q_k$-name for an ordinal. Then there are $n_0$ and a nice norm
$H$ on $n_0$ such that $H(n_0)\geq m$ and 
\begin{quotation}
\noindent if $\t'\geq S_H(\t_0^\uh,\ldots,\t_{n_0-1}^\uh)$, $\|\t'\|>0$,
$w\subseteq w'\subseteq L_0(\langle\rangle)$ 

\noindent then there exist $v\subseteq\cont(\t')$ and $\langle\t_0',\t_1',
\ldots\rangle$ such that 
\[\langle \emptyset,\t_0',\t_1',\ldots\rangle\geq\langle\emptyset,\t_{n_0},
\t_{n_0+1},\ldots\rangle\]
and $\langle w'\cup v,\t_0',\t_1',\ldots\rangle$ decides the value of $\tau$.
\end{quotation}
\end{claim}

\noindent {\em Proof of Claim:}\ \ \ This is essentially 2.14 of
\cite{Sh:207}. 

Define the function $H:\fsuo\lra\omega$ by:  
\begin{description}
\item[$H(u)\geq 0$] always
\item[$H(u)\geq 1$] if $|u|>1$ and for each $\t_i'\geq\t_i^\uh$,
$\|\t^\prime_i\|>0$ (for $i\in u$), for every $w'$, $w\subseteq w'\subseteq
L_0(\langle\rangle)$ there is $v\subseteq\bigcup_{i\in u} \cont(\t_i')$ such
that some pure extension of $\langle w'\cup v,\t_l,\t_{l+1},\ldots\rangle$
decides the value of $\tau$ (where $l=\max u +1$).
\item[$H(u)\geq n+1$] if for every $u'\subseteq u$ either $H(u')\geq n$ or
$H(u\setminus u')\geq n$ (for $n>0$).
\end{description}
As $H$ is monotonic it is enough to find $u$ such that $H(u)\geq m$. The
existence of the $u$ can be proved by induction on $m$, for all
sequences $\langle w,\t_0,\t_1,\ldots\rangle$. 

\noindent Let us start with the case $m=1$. Suppose that $\t'_i\geq \t^\uh_i$,
$\|\t_i^\prime\| >0$ (for $i\in\omega$). For each $i\in\omega$ choose a
creature $\t^*_i\geq\t_i$ such that
\[\cont(\t^*_i)=\cont(\t^\prime_i)\quad\mbox{ and }\quad
\|\t^*_i\|\geq\frac{1}{2}\|\t_i\|\] 
(possible by the definition of the upper half of a creature). Then for each
$w'\subseteq L_0(\langle\rangle)$:
\[\langle w',\t^*_0,\t^*_1,\t^*_2,\ldots\rangle\in\q_k\]
and thus we find $n(w')\in\omega$, $v(w')\in\fsuo$ such that
\[v(w')\subseteq\bigcup_{n<n(w')}\cont(\t^*_n)=\bigcup_{n<n(w')}
\cont(\t^\prime_n)\quad \mbox{ and}\] 
some pure extension of 
\[\langle w'\cup v(w'),\t^*_{n(w')},\t^*_{n(w')+1},\ldots \rangle\]
(and so of $\langle w'{\cup} v(w'),\t_{n(w')},\t_{n(w')+1},\ldots \rangle$)
decides the value of $\tau$.\\
Let $M(\t'_0,\t'_1,\t'_2,\ldots)$ be the first $M$ such that for every 
$w'\subseteq L_0(\langle\rangle)$ there is
$v\subseteq\bigcup_{n<M}\cont(\t'_n)$ such that
\smallskip 

some pure extension of $\langle w'\cup v,\t_{M},\t_{M+1},\ldots
\rangle$ decides the value of $\tau$.
\smallskip

The space 
\[\{\langle\t'_0,\t_1',\t'_2,\ldots\rangle: (\forall i\in\omega)(\t'_i\geq
\t^\uh_i)\}\]
equipped with the natural (product) topology is compact and the function
\[M:\langle\t'_0,\t_1',\t'_2,\ldots\rangle\lra M(\t'_0,\t_1',\t'_2,\ldots)\]
is continuous. Hence the function $M$ is bounded, say by $n_0$. Clearly
$H(n_0)\geq 1$.

\noindent Now suppose that we always can find a set of the norm $\geq m\geq
1$. Thus we find an increasing sequence $\langle l_i: i\in\omega\rangle$ such
that $H([l_i,l_{i+1}))\geq m$ for each $i$. Consider the space of all
increasing $\psi\in\baire$ such that $(\forall i\in\omega)(\psi(i)\in
[l_i,l_{i+1}))$ - it is a compact space. For each $\psi$ from the space we may
consider $\langle w,\t_{\psi(0)},\t_{\psi(1)},\ldots\rangle$ and the
respective function $H_\psi$. By the induction hypothesis we find $n=n_\psi$
such that $H_\psi(n)\geq m$. But $H_\psi(n)\leq H(\{\psi(i): i<n\})$. By the
compactness we find one $n$ such that $m\leq H(\{\psi(i): i<n\})$ for each
$\psi$. Hence we conclude that $H(l_{n+1})\geq m+1$. 

\begin{mainclaim}
\label{cl10}
\begin{enumerate}
\item The forcing notion $\q_k$ has the $S_k$-localization property. 
\item Moreover the following stronger condition is satisfied by $\q_k$:
\begin{description}
\item[$(S^*_k)_{\q_k}$] \ \ \ Suppose that $N$ is a countable
elementary submodel of $\langle {\cal H}(\beth_7^+),\in,<^*\rangle$, 
$p=\langle w,\t_0,\t_1,\ldots\rangle\in\q_k\cap N$. Assume that $Y\in\iso$
is such that 
\[(\forall X\in\iso\cap N)((X,Y)\in S_k).\]
Then there is a condition $q\geq p$ which is $(N,\q_k)$-generic and such
that
\[q\forces_{\q_k}(\forall X\in N[\dot{G}_{\q_k}]\cap\iso)((X,Y)\in S_k).\] 
\end{description}
\end{enumerate}
\end{mainclaim}

\noindent{\em Proof of Main Claim:}\ \ \ 1) Since in the next section we will
need the property $(S^*_k)_{\q_k}$, we will present the proof of it fully
below. Here we sketch the proof of the $S_k$-localization property for readers
not interested in the stronger property (needed for iterations).  So suppose
that $\dot{X}$ is a $\q_k$-name for an element of $\iso$ and $p=\langle
w,\t_0,\t_1,\ldots \rangle\in\q_k$. We may assume that $\|\t_n\|\geq n$. We
inductively define integers 
\[0=b_0<b_1<b_2<\ldots\quad\mbox{ and }\quad 0=n_0<n_1<n_2<\ldots\]
and nice norms $H_m$ on $[n_m,n_{m+1})$. Suppose we have defined $b_m,n_m$.

Let $\tau_m$ be a $\q_k$-name for an integer such that
\[\forces_{\q_k} |[b_m,\tau_m)\cap \dot{X}|>2.\]
Modifying the tail (above $n_m$) of $p$ we may assume that if $n_m\leq n$,
$w\subseteq w'\subseteq R_n(\langle\rangle)$ and some pure extension of
$\langle w',\t_{n+1},\t_{n+2},\ldots\rangle$ decides the value of $\tau_m$
then $\langle w',\t_{n+1},\t_{n+2},\ldots\rangle$  does it already (see
\ref{cl6}). Applying \ref{cl9} we find $n_{m+1}>n_m$ and a nice norm $H_m$ on
$[n_m,n_{m+1})$ such that
\begin{quotation}
\noindent $H_m([n_m,n_{m+1}))\geq m+1$ and if $\t'\geq
S_{H_m}(\t^\uh_{n_m},\ldots,\t^\uh_{n_{m+1}-1})$, $w\subseteq w'\subseteq
L_{n_m}(\langle\rangle)$ then there is $v\subseteq\cont(\t')$ such that some
pure extension of 
\[\langle w'\cup v,\t_{n_{m+1}},\t_{n_{m+1}+1},\ldots\rangle\]
decides the value of $\tau_m$ and thus $\langle w'\cup
v,\t_{n_{m+1}},\t_{n_{m+1}+1},\ldots\rangle$ does it.
\end{quotation}
Let $b_{m+1}$ be an integer larger than all possible values forced to $\tau_m$
in the condition above. 

Now for each $l\in\omega$ we put
\[\t^*_l=S\big(S_{H_{lk}}(\t^\uh_{n_{lk}},\ldots,\t^\uh_{n_{lk+1}-1}),\ldots,
S_{H_{(l+1)k-1}}(\t^\uh_{n_{(l+1)k-1}},\ldots,\t^\uh_{n_{(l+1)k}-1})\big),\]
and then 
\[q=\langle w,\t^*_0,\t^*_1,\ldots\rangle,\quad\quad B=\{b_0,b_1,b_2,
\ldots\}.\]  
Check that $q\geq p$ and
\[q\forces_{\q_k}\mbox{``the set }B\ \ \ S_k\mbox{--localizes
}\dot{X}\mbox{''}\]
(or see the end of the proof of 2) below).
\medskip

\noindent 2)\ \ \ The construction of the condition $q$ required in
$(S^*_k)_{\q_k}$ is similar to that in 1). Here, however, we have to take care
of all names for elements of $\iso$ from the model $N$ (as well as names for
ordinals -- to ensure the genericity).

\noindent Let $\langle\sigma_n: n\in\omega\rangle$ enumerate all $\q_k$-names
from $N$ for ordinals and let $\langle\dot{A}_n: n\in\omega\rangle$ list all
names (from $N$) for infinite subsets of $\omega$. Of course, both sequences
are {\em not} in $N$ but all their initial (finite) segments are there.

Now we inductively define sets $B_n=\{b^n_0,b^n_1,b^n_2,\ldots\}\in\iso$ and
conditions $q_n=\langle w,\t^n_0,\t^n_1,\t^n_0,\ldots\rangle\in\q_k$ such
that $B_n,q_n\in N$:

\noindent To start with we put $B_0=\omega$, $q_0=p=\langle
w,\t_0,\t_1,\ldots\rangle$. 

\noindent Arriving at stage $n>0$ we have defined $B_{n-1}$, $q_{n-1}\in N$.
We define $B_n$, $q_n$ applying the following procedure {\em inside} the
model $N$ (so the result will be there; compare this procedure with that in
part (1)): 
\medskip

\noindent Let $\tau_0$ be a $\q_k$-name for an integer such that
$\forces_{\q_k}(\forall i<n)(|\tau_0\cap\dot{A}_i|>2)$. We modify
$q_{n-1}$ (passing to a pure extension of it) and we assume that 
\begin{description}
\item[$(*)_{\tau_0,\sigma_0,\ldots,\sigma_{n-1}}$]\ \ \ if $v\subseteq
R^{n-1}_i(\langle\rangle)$, $i\in\omega$ and there exists a pure extension of
$\langle w\cup v,\t_{i+1}^{n-1},\t_{i+2}^{n-1},\ldots\rangle$ deciding the
value of one of $\tau_0$, $\sigma_0,\ldots,\sigma_{n-1}$ then $\langle w\cup
v,\t_{i+1}^{n-1},\t_{i+2}^{n-1},\ldots\rangle$ decides it already
\end{description}
(see \ref{cl6}). Next, by Claim \ref{cl9}, we find $n_0$ and a nice norm $H_0$
on $n_0$ such that 
\begin{quote}
\noindent $H_0(n_0)\geq 1$ and if $\t'\geq
S_{H_0}((\t_0^{n-1})^\uh,\ldots,(\t_{n_0-1}^{n-1})^\uh)$, $\|\t'\|>0$, 

\noindent $w\subseteq w'\subseteq L_0^{n-1}(\langle\rangle)$ then there
exists $v\subseteq\cont(\t')$ such that 

\noindent $\langle w'\cup v,\t_{n_0}^{n-1},\t_{n_0+1}^{n-1},\t_{n_0+2}^{n-1} e
\ldots\rangle$ decides the value of $\tau_0$. 
\end{quote}
Let $\t_0^n=S_{H_0}((\t_0^{n-1})^\uh,\ldots,(\t_{n_0-1}^{n-1})^\uh)$ and let
$b_0^n$ be greater than all possible values of $\tau_0$ (i.e. the values
forced in the condition on $H_0$ above). Let $\tau_1$ be a $\q_k$-name for an
integer such that $\forces_{\q_k}(\forall i<n)(|\:[b_0^n,\tau_1)\cap
\dot{A}_i\:|>2)$. We modify ``the tail'' of $q_{n-1}$ and we assume
$(*)_{\tau_1,\sigma_0,\ldots, \sigma_{n-1}}$ (for $i\geq n_0$, $q_{n-1}$).
Next we choose $n_1>n_0$ and a nice norm $H_1$ on $[n_0,n_1)$ such that
\begin{quote}
\noindent $H_1([n_0,n_1))\geq 2$ and if $\t'\geq
S_{H_1}((\t_{n_0}^{n-1})^\uh,\ldots,(\t_{n_1-1}^{n-1})^\uh)$, $\|\t'\|>0$,
$w\subseteq w'\subseteq L_{n_0}^{n-1}(\langle\rangle)$ then there exists
$v\subseteq\cont(\t')$ such that $\langle w'\cup v,\t_{n_1}^{n-1},
\t_{n_1+1}^{n-1},\ldots\rangle$ decides the value of $\tau_1$. 
\end{quote}
Let $\t_1^n=S_{H_1}((\t_{n_0}^{n-1})^\uh,(\ldots,\t_{n_1-1}^{n-1})^\uh)$ and
let $b_1^n$ be greater than all possibilities for $\tau_1$ in the above
property.  

\noindent We continue in this fashion and we determine integers
$n_0<n_1<n_2<\ldots$, $b_0^n<b_1^n<b_2^n<\ldots$ and nice norms $H_0$, $H_1$,
$H_2\ldots$ and we define creatures $\t_i^n=S_{H_i}((\t_{n_{i-1}}^{n-1})^\uh,
\ldots, ((\t_{n_i-1}^{n-1})^\uh)$. Finally we let  
$B_n=\{b_0^n,b_1^n,b_2^n,\ldots\}$ and $q_n=\langle w,\t^n_0,\t^n_1,\ldots
\rangle$ (actually we should have taken more care while getting
$(*)_{\tau_i,\sigma_0,\ldots,\sigma_{n-1}}$ in constructing $\t^n_i$: this is
necessary for $\lim\limits_{i\rightarrow\infty}\|\t^n_i\|=\infty$). 
\medskip

Suppose now that $Y\in\iso$ is a set $S_k$-localizing $\iso\cap N$. Then for
each $n\in\omega$ we have
\[(\exists^\infty i)(\forall m<k)(|\:[\mu_Y(i+m),\mu_Y(i+m+1))\cap B_n|\geq
2).\] 
We inductively define increasing sequences $\langle i_n: n\in\omega\rangle$,
$\langle l_n: n\in\omega\rangle$ of integers and a sequence $\langle \t^*_n:
n\in\omega\rangle$ of creatures:
\medskip

\noindent Let $i_0$ be such that 
\[(\forall m<k)(|\:[\mu_Y(i_0+m),\mu_Y(i_0+m+1))\cap B_1|\geq 2).\]
Let $j^0_m\in\omega$ (for $m<k$) be such that $\mu_{B_1}(j^0_m)$ is the
first element of 
\[[\mu_Y(i_0+m),\mu_Y(i_0+m+1))\cap B_1.\]
Let $\t^*_0=S(\t_{j^0_0+1}^1,\t_{j^0_1+1}^1,\ldots,\t_{j^0_{k-1}+1}^1)$. 
Note that the creatures $\t_j^1$ were obtained as results of the operation
$S_{H_j}$ (for some norms $H_j$) and hence their roots (i.e. $\langle
\rangle$) are $>k$-splitting points, thus no danger can appear in this
procedure. Finally we choose $l_0$ such that
$\t^*_0\in\Sigma^*(\t_0,\ldots,\t_{l_0-1})$. 
\smallskip

\noindent Assuming that we have defined $i_n$, $l_n$, $\t^*_n$, take
$i_{n+1}>i_n+k$ such that
\[(\forall m<k)(|\:[\mu_Y(i_{n+1}+m),\mu_Y(i_{n+1}+m+1))\cap B_{n+2}|\geq
2)\] 
and if $\mu_{B_{n+2}}(j)\in [\mu_Y(i_{n+1}),\mu_Y(i_{n+1}+1))$ then
$\t^{n+2}_j\in\Sigma^*(\t_{l_n},\ldots,\t_m)$ (for some $m>l_n$) and
$\|\t^{n+2}_j\|\geq n+1$. Now, let $j^{n+1}_m$ (for $m<k$) be such that
$\mu_{B_{n+2}}(j^{n+1}_m)$ is the first element of $[\mu_Y(i_{n+1}+m),
\mu_Y(i_{n+1}+m+1))\cap B_{n+2}$. Finally,
$\t^*_{n+1}=S(\t^{n+2}_{j_0^{n+1}+1},\ldots, \t^{n+2}_{j^{n+1}_{k-1}+1})$ 
and $l_{n+1}$ is such that
\[\t^*_{n+1}\in\Sigma^*(\t_{l_n},\ldots,\t_{l_{n+1}-1}).\]  
The condition $q$ is $\langle w,\t^*_0,\t^*_1,\ldots\rangle$.
Clearly $q\geq p$. To show that it is $(N,\q_k)$-generic suppose that
$\sigma\in N$ is a $\q_k$-name for an ordinal, say $\sigma=\sigma_n$, and let
$q'=\langle w',\t_0',\t_1',\ldots\rangle\geq q$ be a condition deciding
$\sigma$. Look at the construction of $q_{n+1}=\langle
w,\t^{n+1}_0,\t^{n+1}_1,\ldots \rangle$. Because of
$(*)_{\tau_0,\sigma_0,\ldots,\sigma_n}$, if $v\subseteq
M^{n+1}_i(\langle\rangle)$ and there is (in $N$) a pure extension of 
$\langle w\cup v,\t^{n+1}_{i+1},\t^{n+1}_{i+2},\ldots\rangle$ deciding
$\sigma_n$ then $\langle w\cup v,\t^{n+1}_{i+1},\t^{n+1}_{i+2},
\ldots\rangle$ decides it already.   

\noindent Let $i$ be such that $w'\subseteq M^{n+1}_i(\langle\rangle)$. Then
there exists a pure extension of the condition $\langle
w',\t^{n+1}_{i+1},\t^{n+1}_{i+2},\ldots\rangle$ deciding $\sigma_n$, e.g.
that one which can be obtained from $q'$ (note that $q'$ and $\langle
w',\t^{n+1}_{i+1},\t^{n+1}_{i+2},\ldots\rangle$ are compatible). By the 
elementarity of $N$ there is such an extension in $N$. This implies that
$\langle w',\t^{n+1}_{i+1},\t^{n+1}_{i+2},\ldots\rangle$ decides $\sigma_n$
(and ``the decision'' belongs to $N$). The conditions $q'$ and $\langle w',
\t^{n+1}_{i+1},\t^{n+1}_{i+2},\ldots\rangle$ are compatible, so the values
given by them to $\sigma_n$ are the same and we are done (with the
genericity). 

Now we want to show that
\[q\forces_{\q_k}(\forall X\in N[\dot{G}_{\q_k}]\cap\iso)((X,Y)\in S_k).\]
Let $\dot{A}\in N$ be a $\q_k$-name for an infinite subset of $\omega$, say
$\dot{A}=\dot{A}_m$. We are going to prove that 
\[q\forces_{\q_k}(\exists^\infty n\in\omega)(\forall i<k)(|\;[\mu_Y(n+i),
\mu_Y(n+i+1))\cap\dot{A}_m\;|\geq 2).\]
Let $q'=\langle w',\t_0',\t_1',\ldots\rangle\geq q$. Fix $l>m$. Look at
$\t_l'$ -- for some $m'$, $m''$, $l<m'<m''$ we have
$\t_l'\in\Sigma^*(\t_{m'}^*,\ldots,\t_{m''}^*)$ and
$w'\subseteq\bigcup_{i<m'}\cont(\t^*_i)$.  So we find $n>l$ and $t\in\t_l'$
such that $\t_l'$ ``above $t$'' comes from $\t^*_n$ by decreasing norms (like
in the operation of taking the upper half but possibly with other values
subtracted) and refining. Note that neccessarily
$|\suc_{\t_l'}(t)|=k$ and $\t'_l$ above each successor of $t$ refines some
$\t_{j^n_i+1}^{n+1}$ (for $i<k$) modulo decreasing norms by some values, which
does not influence contributions. Next we find $v\subseteq\cont(\t_l')$
(actually $v=v_0\cup\ldots\cup v_{k-1}$, $v_i$ is included in the
contributions of the part above that successor of $t$ which refines
$\t_{j^n_i+1}^n$) such that $\langle w'\cup v,\t_{l+1}',\t_{l+2}',\ldots
\rangle$ decides the values of the names $\tau_{j^n_0+1},\ldots,
\tau_{j^n_{k-1}+1}$ (defined in the construction of $\t^{n+1}_j$'s) and forces
them to be bounded by $b_{j^{n+1}_0+1},\ldots,b_{j^{n+1}_{k-1}+1}$,
respectively. By the choice of $\tau_j$'s this means that
\[\langle w'\cup v,\t_{l+1}',\t_{l+2}',\ldots\rangle\forces_{\q_k}(\forall
i<k)(\forall j<n)(|\:[b_{j^{n+1}_i}^n,b_{j^{n+1}_i+1})\cap\dot{A}_j\:|\geq
2).\]  
Since $m<l\leq n$ and $\mu_Y(i_n+i)\leq\mu_{B_{n+1}}(j^n_i)<
\mu_{B_{n+1}}(j^n_i+1) <\mu_Y(i_n+i+1)$ (for all $i<k$) we get 
\[\langle w'\cup v,\t_{l+1}',\t_{l+2}',\ldots\rangle\forces_{\q_k}(\forall
i<k)(|\:[\mu_Y(i_n+i),\mu_Y(i_n+i+1))\cap\dot{A}_m\:|\geq 2).\]   
As $l>m$ was arbitrary we are done. The Main Claim and the theorem are
proved.\QED 
\medskip

\noindent {\sc Remark:} More examples of forcing notions distinguishing the
localization properties we have introduced in this section will be presented
in \cite{RoSh:470}. They are (like the forcing notion $\q_k$) applications of
the general schemata of that paper.
\medskip

\section{Cardinal coefficients related to the localizations}

In this section we discuss cardinal coefficients related to the
localization properties introduced earlier.

Following Vojt\'a\v s (cf \cite{Voj}) with any relation $R\subseteq
X\times Y$ we may associate two cardinal numbers (the unbounded and the
dominating number for $R$):
\[\unbounded(R)=\min\{|B|:(\forall y\in Y)(\exists x\in B)((x,y)\notin R)\}\]
\[\dominating(R)=\min\{|D|:(\forall x\in X)(\exists y\in D)((x,y)\in R)\}.\]
For purposes of applications these cardinals are introduced for
relations $R\subseteq X\times Y$ such that $\dom(R)=\dom(cR)=X$
and $\rng(R)=\rng(cR)=Y$. Note that for each such relation we
have
\[\unbounded(R)=\dominating(cR^{-1}),\ \ \dominating(R)=\unbounded(cR^{-1}),\
\ \unbounded(R^{-1})=\dominating(cR),\ \
\dominating(R^{-1})=\unbounded(cR).\] 
All results of the previous sections provide information on dominating
numbers $\dominating(R)$ for the considered relations. Let
\[\unbounded=\unbounded(\leq^*)=\min\{|F|: F\subseteq\baire\ \&\ (\forall
x\in\baire)(\exists y\in F)(\exists^\infty n\in\omega)(x(n)<y(n))\}\]
\[\dominating=\dominating(\leq^*)=\min\{|F|: F\subseteq\baire\ \&\ (\forall
x\in\baire)(\exists y\in F)(\forall^\infty n\in\omega)(x(n)\leq y(n))\}\]
\[\non(\k)=\min\{|X|: X\subseteq\can\ \&\ X \mbox{ is not meager }\}.\]
These are three of ten cardinal invariants forming the Cicho\'n
Diagram. For more information on the cardinals related to measure and
category see \cite{BJSh:368} or \cite{CiPa}. 

\begin{corollary}\
\begin{description}
\item[1. (see \ref{trivial})\ \ \ ] $(\forall
k\in\omega)(\dominating(\rall)=\con)$ 
\item[2. (see \ref{unbounded})\ \ \ ] $(\forall k>0)(\dominating(c\rall) =
\unbounded)$ 
\item[3. (see \ref{dominujaca})\ \ \ ] $(\forall
k\in\omega)(\dominating(c\rexists)=\dominating)$
\item[4. (see \ref{exists})\ \ \ ] $(\forall
k\in\omega)(\dominating(\rexists) =\dominating(R^\exists_0))$
\item[5. (see \ref{meabou})\ \ \ ] $\unbounded\leq\dominating(R^\exists_0)
\leq\non(\k)$\QED 
\end{description}
\end{corollary}

\begin{corollary}\
\begin{description}
\item[1.] {\sc Con}($\dominating(R^\exists_0)<\dominating$)
\item[2.] {\sc Con}($\dominating<\dominating(R^\exists_0)$)
\end{description}
\end{corollary}

\Proof For the first model add $\aleph_2$ Cohen reals to a model of CH. Then
in the extension we will have $\dominating=\aleph_2$ and
$\dominating(R^\exists_0)=\aleph_1$ (see \ref{cora}). The second model can
be obtained by adding $\aleph_2$ random reals to a model of CH, which results
in a model for $\dominating=\aleph_1$ and
$\dominating(R^\exists_0)=\aleph_2$ (see \ref{cora}). \QED

\begin{corollary} Let $\phi\in\baire$ be an increasing function. Then:
\begin{description}
\item[1. (see \ref{ogolne})\ \ \ ] $\unbounded=\dominating(S_1)\leq
\dominating(S_2) \leq\ldots\leq \dominating(S_k)\leq\ldots\leq
\dominating(S_+)$,\\
$\dominating(S_+)\leq\dominating(S_{+\epsilon}), \dominating(S_+^\phi)
\leq\dominating$, and\\ 
if $\phi$ is increasing fast enough (e.g. $\phi(n)>2^{2n}$) then
$\dominating(S_{+\epsilon})\leq\dominating(S^\phi_+)$,
\item[2. (see \ref{non})\ \ \ ] $\dominating(S^\phi_+)\leq\non(\k)$. \QED. 
\end{description}
\end{corollary}

\begin{corollary}
Let $\phi\in\baire$ be an increasing function. Then
\[\mbox{{\sc Con}}(\dominating(S^\phi_+)=\aleph_1 + \non(\k)=\dominating
=\aleph_2).\]
\end{corollary}

\Proof Start with a model for CH and add to it first $\aleph_2$ Cohen
reals (what causes $\dominating=\aleph_2$ but keeps
$\dominating(S^\phi_+)=\aleph_1$) and next add $\aleph_2$ random reals (what
preserves $\dominating=\aleph_2$, $\dominating(S^\phi_+)=\aleph_1$ but
causes $\non(\k)=\aleph_2$), see \ref{cos}. \QED

\begin{corollary}
Let $k>0$. Then
\[\mbox{{\sc Con}}(\dominating(S_k)<\dominating(S_{k+1})).\]
\end{corollary}

\Proof Let $\q_k$ be the forcing notion from the proof of \ref{eska}. It is
proper (see \ref{cl6}). To get the respective model it is enough to take the
countable support iteration of the length $\omega_2$ of the forcing notions
$\q_k$ over a model of CH. As $\q_k$ does not have the $S_{k+1}$-localization
we easily get that in the resulting model we will have
$\dominating(S_{k+1})=\aleph_2$. The only problem is to show that the
iteration has the $S_k$-localization property (to conclude that in the 
extension $\dominating(S_k)=\aleph_1$). But this is an application of \S3
Chapter XVIII of \cite{Sh:f}. We may think of $S_k$ as a relation on $\baire$
(after canonical mapping). Keeping the notation of \cite{Sh:f} we put:
\medskip

\noindent $S\subseteq{\cal S}_{<\aleph_1}(H(\aleph_1)^{{\bf V}_1})$ and for
$a\in S$, ${\bf g}_a\in\baire$ is such that for each $f\in a\cap\baire$,
$f\; S_k\; {\bf g}_a$,\ \  ${\bf g}=\langle {\bf g}_a: a\in S\rangle$;\\ 
$\alpha^*=1$; $\bar{R}=\langle R_0\rangle=\langle S_k\rangle$.
\medskip

\noindent Note that $(\bar{R},S,{\bf g})$ strongly covers iff it covers iff
$S$ is stationary (for ``strongly covers'' we are in Possibility B). The
property $(S^*_k)_{\q_k}$ of claim~\ref{cl10}(2) guarantees that the forcing
notion $\q_k$ is $(\bar{R},S,{\bf g})$-preserving. Hence Theorem 3.6 of
Chapter XVIII of \cite{Sh:f} applies to this situation and the iteration is
$(\bar{R},S,{\bf g})$-preserving. Consequently we are done. \QED
\smallskip

\begin{problem}
Are the following consistent:
\begin{enumerate}
\item $\unbounded<\dominating(S_2) <\dominating(S_3)< \ldots\dominating(S_k)
<\dominating(S_{k+1}) <\ldots$ ?
\item There exists a sequence
$\langle\phi_\alpha:\alpha<\omega_1\rangle\subseteq\baire$ of increasing 
functions such that 
\[\alpha<\beta<\omega_1\ \Rightarrow\
\dominating(S_+^{\phi_\alpha})<\dominating(S_{+}^{\phi_\beta})?\] 
\end{enumerate}
\end{problem}

\end{document}